\documentclass[11pt]{amsart} 
\usepackage{setspace} \onehalfspace

\usepackage[all]{xy} \usepackage{amsmath} \usepackage{amssymb}
\usepackage{amsthm} 
\def\cal{\mathcal}

\newcommand{\trans}[1]{\hspace{0,1cm}
  #1^{\hspace{-0,4cm}t}\hspace{0,24cm}}

\newcommand{\R}{\mathbb{R}}
\newcommand{\Z}{\mathbb{Z}}

\def\d{\mbox{d}}

\def\Pic{\mbox{Pic}}

\def\Tors{\mbox{Tors}}

\begin{document}

\title{A note on the cone of mobile curves}
\author{Matei Toma} \date{\today}
\thanks{ AMS
  Classification (2000): 32J25; Secondary: 32J13, 14C30.}
\address{Institut de Math\'ematiques Elie Cartan,
Nancy-Universit\'e,\\
B.P. 239, 54506 Vandoeuvre-l\`es-Nancy Cedex, France}
% \\ and Institute of Mathematics of the Romanian Academy.}
%\setcounter{section}{0}
%\setcounter{tocdepth}{1}
\maketitle

\begin{abstract}
S. Boucksom, J.-P. Demailly, M. P\u aun and Th. Peternell proved that the cone of mobile curves $\overline{ME(X)}$ of a projective complex manifold $X$ is dual to the cone generated by classes of effective divisors and conjectured an extension of this duality in the K\"ahler set-up. We show that their conjecture implies that $\overline{ME(X)}$ coincides with the cone of integer classes represented by closed positive smooth $(n-1,n-1)$-forms. Without assuming the validity of the conjecture we prove that this equality of cones still holds at the level of degree functions.
\end{abstract}

\noindent

Let $X$ be a smooth complex projective variety of dimension $n$. A curve $C$ on $X$ is called {\em mobile} if it is the member of an algebraic family of (generically) irreducible curves covering $X$. Let $\overline{ME(X)}$ denote the closed convex cone generated by classes of mobile curves inside 
$N_1(X):=(H^{n-1,n-1}_{\R}(X)\cap H^{2n-2}(X,\Z)/\Tors)\otimes_{\Z}\R$. We shall call the elements of $\overline{ME(X)}$
{\em mobile classes}. 

In \cite{BDPP} it is shown that the following cones in $N_1(X)$ coincide:
\begin{enumerate}
 \item
 the cone  $\overline{ME(X)}$ of mobile curves, 
\item 
the cone $\cal M_{NS}:=\cal M\cap N_1$,\\
where $\cal M\subset H^{n-1,n-1}_{\R}(X)$ is the closure of the convex cone generated by cohomology classes of currents of the type
$\nu_*(\tilde\omega_1\wedge...\wedge\tilde\omega_{n-1})$ for K\"ahler forms  $\tilde\omega_1$, ..., $\tilde\omega_{n-1}$ on a modification $\nu:\tilde X\to X$ de $X$,
%\item
%the closed cone $\overline{SME(X)}$ of strongly mobile curves,\\
%where $SME(X)$ is the convex cone genearted by curves $C=\nu_*(\tilde A_1\cap...\cap\tilde A_{n-1})$ which appear as images of complete intersections of ample divisors $\tilde A_1$, ..., $\tilde A_{n-1}$ on modifications $\nu:\tilde X\to X$,
%\item
%the closed cone  $\overline{ME_{nef}(X)}$ of curves with nef normal bundle,
\item
the dual cone $(\cal E_{NS})^{\vee}$  of the cone    $\cal E_{NS}$ of pseudo-effective divisors on $X$.
\end{enumerate} 

It was known from \cite{Dem92} that  $\cal E_{NS}=  \cal E\cap NS_{\R}$, where   $\cal E$  is the cone of classes of positive closed currents of type $(1,1)$ and $NS_{\R}(X):=(H^{1,1}_{\R}(X)\cap H^{2}(X,\Z)/\Tors)\otimes_{\Z}\R.$
In \cite{BDPP} Conj. 2.3 it is further conjectured that the cones $\cal M$ and  $\cal E$ are dual. 

In this note we compare $\overline{ME(X)}$ to the closed convex cone $P^{n-1,n-1}$ in  $H^{n-1,n-1}_{\R}(X)$ generated by closed positive smooth $(n-1,n-1)$-forms. From the above statements it is clear that $\overline{P^{n-1,n-1}}\cap N_1(X)\subset \overline{ME(X)}$. The converse inclusion will follow from our arguments if we admit the conjecture of  Boucksom,  Demailly,  P\u aun,  Peternell. If not we still get an equality at the level of degree functions as follows.

Any mobile class $\alpha$ gives rise to a degree function  
$$\deg_{\alpha}: \Pic(X)\to\R, \ \ \deg_{\alpha}L:=c_1(L)\alpha$$
and further to a notion of stability of torsion-free sheaves on $X$ which generalizes the classical case 
$\alpha=H^{n-1}$ for a class $H$ of an ample divisor, cf. \cite{CaPe}. On the other side we consider {\em semi-K\"ahler} metrics on $X$, i.e. such that their associated K\"ahler forms $\omega$ satisfy $\d\omega^{n-1}=0$. Such a metric %(in fact any hermitian metric)
 gives likewise a degree function:
$$\deg_{\omega}: \Pic(X)\to\R, \ \ \deg_{\omega}L:=c_1(L)[\omega^{n-1}].$$
Then we can state:

{\bf Theorem.}
{\it For any mobile class $\alpha$ in the interior of the mobile cone $\overline{ME(X)}$ there exists a semi-K\"ahler metric on $X$ with associated K\"ahler form $\omega$ such that
$\deg_{\alpha}=\deg_{\omega}$ on $\Pic(X)$.}

By \cite{LT} Cor. 5.3.9 this has the following consequence on the moduli space of stable vector bundles:

{\bf Corollary.}
{\it If $\alpha$ is a class in the interior of $\overline{ME(X)}$, 
then the smooth part of the moduli space of stable vector bundles with respect to $\alpha$ admits a natural K\"ahler structure.}

The following notation will be used:\\
We denote by
 $\cal D'^{p,q}$ the space of currents of bidegree $(p,q)$ (and bidimension $(n-p,n-q)$) on $X$. Let 
$$V_{\d}=V_{\d}(X):=\{ T\in \cal D'^{1,1}_{\R}\ | \ \d T=0\}/\d\d^c\cal D'^{0,0}_{\R},$$
$$V_{\d\d^c}=V_{\d\d^c}(X):=\{ T\in \cal D'^{1,1}_{\R}\ | \ \d\d^c T=0\}/\{ \overline\partial S+\partial \bar S\ | \
S\in\cal D'^{1,0}\},$$
$V_{\d}^+$, $V_{\d\d^c}^+$ the cones generated by positive currents in  $V_{\d}$ and $V_{\d\d^c}$ respectively and \\
$V_{\d, NS}^+$, $V_{\d\d^c, NS}^+$ their intersections with
$NS_{\R}(X):=(H^{1,1}_{\R}(X)\cap H^{2}(X,\Z)/\Tors)\otimes_{\Z}\R.$

We first prove a

{\bf Lemma.} {\it
The natural map $j:V_{\d}\to V_{\d\d^c}$ induces a bijection between  positive cones
$$V_{\d, NS}^+\to V_{\d\d^c, NS}^+.$$}

\begin{proof}
It is known that in the case of compact K\"ahler manifolds the map $j$, which associates to a class of a bidegree $(1,1)$
closed current $\{T\}_{\d}\in V_{\d}$ its class $\{T\}_{\d\d^c}\in V_{\d\d^c}$, is well defined and bijective.

The inclusion
$j(V_{\d}^+)\subset V_{\d\d^c}^+$  being obvious, we consider a $\d\d^c$-closed positive current  $T$ with 
  $\{T\}_{\d\d^c}\in V_{\d\d^c, NS}$. Let $\eta:=j^{-1}( \{T\}_{\d\d^c})\in V_{\d}$. We shall show that $\eta\in V_{\d}^+$.
  
By the cited result of \cite{BDPP} it suffices to check that for any 
  modification
$\nu:\tilde X\to X$ and K\"ahler forms $\tilde\omega_1$, ...,$\tilde\omega_{n-1}$ on $\tilde X$ one has
$$\eta\nu_*([\tilde\omega_1\wedge...\wedge\tilde\omega_{n-1}])\ge 0.$$

But
 $$\eta\nu_*([\tilde\omega_1\wedge...\wedge\tilde\omega_{n-1}])=\nu^*(\eta)[\tilde\omega_1\wedge...\wedge\tilde\omega_{n-1}]=
\nu^*(\{T\}_{\d\d^c})[\tilde\omega_1\wedge...\wedge\tilde\omega_{n-1}]$$
and Theorem 3
of \cite{AlBa95} asserts that
$\nu^* (\{T\}_{\d\d^c})\in V_{\d\d^c}^+(\tilde X)$, whence the desired inequality.
\end{proof}
%%%%%%%%%%%%%%%%%%%%%%%%%%%%%%%%%%%%%%%%%%%%%%%%%%%%%%%%%%%%%%%%%%%%%%%%%%%%%%%%%%%%%%%%%%%%%%%%%%%%%%%%%%%%%%%%%%%%%%%%%%

We can now give the proof of the Theorem.
\begin{proof}
Let $\alpha$ be an element in the interior of the cone of mobile curves.
We denote by $\cal D'^{1,1}_+$ the cone of positive currents inside the space
 $\cal D'^{1,1}$ of bidegree $(1,1)$ currents on $X$. We fix a K\"ahler form  
$\sigma$  on $X$ and set $\cal D'^{1,1}_{+,\sigma}:=\{ T\in \cal D'^{1,1}_+ \ | \ \int_X T\wedge\sigma^{n-1}=1\}$.
This is a  compact set for the weak topology on  $\cal D'^{1,1}$,  \cite{D} III.1.23. 
Let $\beta_1,...,\beta_k\in H^{n-1,n-1}_{\R}(X)$ be such that 
$V_{\d\d^c, NS}=\{t\in V_{\d\d^c} \ | \ t\beta_1=0, ...,t\beta_k=0\}$ and set
$W:=\{ T\in \cal D'^{1,1} \ | \ \d\d^c T=0, \ \{T\}_{\d\d^c}\alpha=0, \ \{T\}_{\d\d^c}\beta_1=0,..., \{T\}_{\d\d^c}\beta_k=0\}$. Remark that $W$ and  $\cal D'^{1,1}_{+,\sigma}$
are disjoint. Indeed, if $T\in \cal D'^{1,1}_{+,\sigma}$ were $\d\d^c$-closed and $\{T\}\in V_{\d\d^c, NS}^+$, then by the Lemma  there would exist a  $\d$-closed positive current
$S\in  \cal D'^{1,1}_+$ such that $\{T\}_{\d\d^c}=j(\{S\}_{\d})$ and in particular
 $$\{T\}_{\d\d^c}\alpha=
\{S\}_{\d}\alpha>0.$$ 
%Donc le sous-espace ferm\'e 
%$W:=\{ T\in \cal D'^{1,1} \ | \ \d\d^c T=0, \ \{T\}_{\d\d^c}\alpha=0\}$ 
%dans $\cal D'^{1,1}$ et $\cal D'^{1,1}_{+,\sigma}$
%sont disjoints. 
The Hahn-Banach theorem then implies the existence of a functional on  $\cal D'^{1,1}$, which vanishes on $W$ and is  positive on $\cal D'^{1,1}_{+,\sigma}$. 
This functional is thus given by a  real $(n-1,n-1)$-form  $u$ on $X$. 
The form $u$ is strictly positive on  $X$ since the functional is positive onr $\cal D'^{1,1}_{+,\sigma}$. The vanishing on $W$ implies that $u$ eis also  $\d$-closed.
As functionals on  $ NS_{\R}(X)$, $[u]$  and $\alpha$ have the same kernel, hence they coincide up to some multiplicative constant.

It is enough now to take a postive  $(n-1)$-st root $\omega$ of $u$. For the convenience of the reader we show how this is done. Remark first that 
     $$(i\sum_{1\le i,j\le n}a_{ij}\d z_i\wedge\d \overline z_j)^{n-1}=
     (n-1)!i^{(n-1)^2}\sum_{1\le i,j\le n}(-1)^{i+j}c_{ji}\hat{\d z_i}\wedge\hat{\d \overline z_j},\leqno{(1)},$$
     where we denoted by
     $c_{ij}$ the cofactor of $a_{ij}$ in the matrix $A=(a_{ij})_{1\le i,j\le n}$ and  
     $\hat{\d z_i}:=\d z_1\wedge...\wedge\d z_{i-1}\wedge\d z_{i+1}\wedge...\wedge\d z_n$, 
     $\hat{\d \overline z_j}
     :=\d \overline z_1\wedge...\wedge\d\overline z_{j-1}\wedge\d\overline z_{j+1}\wedge...\wedge\d \overline z_n$.
     The relation $\trans{C}A=\det (A)I_n$ for the matrix of cofactors $C=(c_{ij})_{1\le i,j\le n}$ implies 
     $$A=\sqrt[n-1]{\det(C)}\trans{C}^{-1}$$
     when  $A$ is positive definite. 
When the matrix  $C$ is positive definite, one obtains an unique positive definite solution
       $A$ of equation (1). 

 \end{proof}

%%%%%%%%%%%%%%%%%%%%%%%%%%%%%%%%%%%%%%%%%%%%%%%%%%%%%%%%%%%%%%%%%%%%%%%%%%%%%%%%%%%%%%%%%%%%%%%%%%%%%%%%%%%%%
%%%%%%%%%%%%%%%%%%%%%%%%%%%%%%%%%%%%%%%%%%%%%%%%%%%%%%%%%%%%%%%%%%%%%%%%%%%%%%%%%%%%%%%%%%%%%%%%%%%%%%%%%%%%%

\end{document}